\documentclass[10pt,oneside]{article}

\usepackage[left=1in,right=1in,top=1in,bottom=1in]{geometry}
\usepackage{authblk}

\usepackage{setspace}
\usepackage[parfill]{parskip}
\usepackage{lscape}

\usepackage{color}
\usepackage[colorlinks,
            citecolor=blue,
            urlcolor=blue,
            linkcolor=blue]{hyperref}

\usepackage[utf8]{inputenc} 
\usepackage[T1]{fontenc}    
\usepackage{hyperref}       
\usepackage{url}            
\usepackage{booktabs}       
\usepackage{amsfonts}       
\usepackage{nicefrac}       
\usepackage{enumitem}
\usepackage{multirow}
\usepackage{algorithm}
\usepackage{nameref}
\usepackage{mathrsfs}
\usepackage{comment}
\usepackage{tablefootnote}
\usepackage[font=footnotesize]{caption}
\usepackage{mathtools}

\usepackage{libertine}
\usepackage[T1]{fontenc}
\usepackage[utf8]{inputenc}

\usepackage{microtype}
\usepackage{graphicx}
\usepackage{subfig}
\usepackage{amsmath}   
\usepackage{amsthm}
\usepackage{amssymb}
\usepackage{siunitx}   
\usepackage{bm}        
\usepackage{thmtools}  
\usepackage{thm-restate}
\usepackage{color}     
\usepackage{xspace}    

\newcommand{\bb}{\bm{b}}
\newcommand{\bc}{\bm{c}}

\newcommand{\bfvec}{\bm{f}}
\newcommand{\br}{\bm{r}}
\newcommand{\bu}{\bm{u}}

\newcommand{\bw}{\bm{w}}
\newcommand{\bx}{\bm{x}}
\newcommand{\by}{\bm{y}}
\newcommand{\bz}{\bm{z}}

\newcommand{\bA}{\bm{A}}

\newcommand{\bD}{\bm{D}}

\newcommand{\bH}{\bm{H}}
\newcommand{\bI}{\bm{I}}

\newcommand{\bL}{\bm{L}}

\newcommand{\bgamma}{\bm{\gamma}}
\newcommand{\bPi}{\bm{\Pi}}
\newcommand{\bphi}{\bm{\phi}}

\newcommand{\E}{\mathbb{E}}

\newcommand{\bzero}{\bm{0}}
\newcommand{\bone}{\bm{1}}

\definecolor{ForestGreen}{RGB}{34,139,34}

\usepackage{pifont}

\newtheorem{proposition}{Proposition}

\usepackage[natbib=true,
            maxcitenames=2,
            maxbibnames=20,
            mincrossrefs=5,
            citestyle=authoryear,
            bibstyle=authoryear,
            uniquename=false,
            uniquelist=false]{biblatex}
\DefineBibliographyStrings{english}{
    andothers = {\mkbibemph{et\addabbrvspace al\adddot}}
}

\renewbibmacro{in:}{}
\DeclareCiteCommand{\citeyearparlink}
    {}
    {\mkbibparens{\bibhyperref{\printdate}}}
    {\multicitedelim}
    {}

\title{Optimal accuracy for linear sets of equations with the graph Laplacian}

\author[1]{Richard B. Lehoucq}
\author[2]{Michael Weylandt}
\author[1]{Jonathan W. Berry}

\affil[1]{Sandia National Laboratories}
\affil[2]{Zicklin School of Business, Baruch College, CUNY}

\usepackage{lineno}

\date{\today}

\begin{document}

\maketitle

\begin{abstract}
We show that certain Graph Laplacian linear sets of equations exhibit optimal accuracy, guaranteeing that the relative error is no larger than the norm of the relative residual and that optimality occurs for carefully chosen right-hand sides. Such sets of equations arise in PageRank and Markov chain theory. We establish new relationships among the PageRank teleportation parameter, the Markov chain discount, and approximations to linear sets of equations. The set of optimally accurate systems can be separated into two groups for an undirected graph---those that achieve optimality asymptotically with the graph size and those that do not---determined by the angle between the right-hand side of the linear system and the vector of all ones.  We provide supporting numerical experiments. 
\end{abstract}
\onehalfspace

\begin{refsection}


An approximate solution $\hat{\bx}$ for a linear set of equations $\bL\bx=\bb$
has roughly $d$ digits of accuracy when $\|\hat{\bx} - \bx\|/\|\bx\| \approx 10^{-d}$.
The classical two-sided inequality
\begin{align}\label{unif-stab}
\frac{1}{\kappa(\bL)} \frac{\|\bb-\bL\hat{\bx} \|}{\|\bb\|} \leqslant \frac{\|\bx - \hat{\bx}\|}{\|\bx\|} \leqslant \kappa(\bL) \frac{\|\bb-\bL\hat{\bx}\|}{\|\bb\|}\quad \bx,\bb \neq 0
\end{align}
implies that the relative error is norm-equivalent to the relative
residual error $\|\bb-\bL\hat{\bx} \|/\|\bb\|$ with constants given by the condition number $\kappa(\bL)$
and its reciprocal. For us, the matrix $\bL$ is a graph Laplacian and the vector $x$ represents a
network centrality measure indicating the importance of the vertices, \emph{e.g.}, the PageRank~\citep{brpa:98}
vector or the vector of mean hitting-times. \
Unfortunately, the condition number $\kappa(\bL)$ increases with  graph size or with the PageRank teleportation parameter rendering \eqref{unif-stab} useless in practice. 
We establish improved variants of the two-sided inequality and explore their profound computational implications.
 
We focus our analysis on the relationship between the relative error and the
relative residual. 
This relationship is key to assessing the quality of $\hat{\bx}$ because it relates the observable quantity
$\|\bL\hat{\bx} -\bb\|$ with the unobservable quantity $\|\hat{\bx} -\bx\|$. 
We show that the strength of this relationship is determined
by the angle between  $\bb$ and the vector of all ones on an undirected graph. 
This relationship is also
dependent on the so-called \emph{Markov chain discount}, a classical concept
that we find is equivalent to the PageRank teleportation parameter for undirected graphs. This
provides an elegant probabilistic basis for the degree-normalized PageRank variant
and the simple characterization of PageRank on an undirected graph sought by \citet[p.356]{glei:15}.
 
Our contributions are twofold: i) we establish a more informative variant of \eqref{unif-stab} using a \emph{data-dependent condition
number} in \S\ref{unif-opt} and ii) we reframe graph centrality measures in the language
of discrete potential theory and show how certain potentials can
achieve asymptotically optimal accuracy in \S\ref{sec:dpt}. 
We discuss the application of these
results to PageRank  in \S\ref{pgrk} and conclude with
numerical simulations highlighting the impact of our improved bounds in \S\ref{experiments}.
 
\section{Graph Condition Number Bounds via Discounting}\label{sec:trans-mat-norm}

Let $\bA$ be an adjacency matrix for a possibly directed, possibly weighted graph, not assumed connected, with $n$ vertices and $\bD$ be the diagonal matrix of row degrees, which we assume to be nonsingular (so that isolated vertices contain a self-loop). The Markov chain transition matrix for the random walk over the graph is then given by $\bD^{-1}\bA$ and is row-stochastic. 
Let $\bL=\bI - \alpha \bD^{-1}\bA$ denote the $\alpha$-discounted graph Laplacian for fixed $0 < \alpha <  1$. Recall the
well-known relationship  $\bD\bone =\bA\bone$, which in turn implies the useful relationship 
\begin{align} \label{constant}
\bL \bone = \bone (1-\alpha) \text{ where } \bone^\top = (1, 1, \cdots, 1)\,.
\end{align} 
This indicates that the angle between $\bb$ and $\bone$ and the angle between $\bx$ and $\bone $ are closely related.
The inequality $\| \bL^{-1} \|_p \leqslant 1/(1-\alpha \|\bD^{-1}\bA\|_p)$ is well-known for the 
matrix $p$-norm \citep[\emph{e.g.}][p.74]{govl:13} and implies that the inverse of $\bL$  is bounded.  
Additionally, when the discount $\alpha$ satisfies $0 < \alpha \|\bD^{-1}\bA\|_p < 1 $, the condition number also satisfies
\begin{align} \label{cond-num}
\kappa_p (\bL) &= \|\bL \|_p \, \| \bL^{-1} \|_p  \leqslant \frac{1+\alpha\|\bD^{-1}\bA\|_p}{1-\alpha\|\bD^{-1}\bA\|_p } \,. 
\end{align}
The condition number is \emph{uniformly} bounded when the previous inequality holds independent of the number of  graph vertices.

A bound on the relative maximum element-wise error occurs with the $p=\infty$ norm. Then $\|\bL\|_\infty=1+\alpha$ and $\|\bD^{-1}\bA\|_\infty =1$ so that $\kappa_\infty (\bL) \leqslant (1+\alpha)/(1-\alpha)$ and
\eqref{unif-stab} implies that
\begin{align} \label{inf-2s}
\frac{1-\alpha}{1+\alpha}\,\frac{\|\br\|_\infty}{\|\bb\|_\infty} \leqslant \frac{\|\bx -\hat{\bx}\|_\infty}{\|\bx\|_\infty} 
\leqslant \frac{1+\alpha}{1-\alpha}\,\frac{\|\br\|_\infty}{\|\bb\|_\infty}
\end{align}
holds for any $0 < \alpha < 1$ and independently of the number of graph vertices where $\|\bx\|_\infty = \max_i |x_i|$ and $\br=\bb-\bL\hat{\bx}$. 

An important case in practice is given by an undirected graph, \emph{i.e.},  $\bA=\bA^\top$ because the change of coordinates $\by=\bD^{1/2}\bx$ on the linear system $\bL\bx=\bb$ leads to the symmetric positive definite system  
\begin{equation}\label{spd}
\bH\by = \bfvec\,, \quad \bH \coloneq \bD^{1/2} \bL \bD^{-1/2}\,,\, \bfvec  \coloneq  \bD^{1/2} \bb\,,
\end{equation}
when $0 <\alpha < 1$.
The bound in \eqref{cond-num} implies that  $\kappa_{\bD} (\bL) =  \kappa_2 (\bH) \leqslant (1+\alpha)/(1-\alpha) $ since 
$\|\bD^{-1}\bA\|_{\bD} = \|\bD^{-1/2}\bA \bD^{-1/2}\|_2 = \rho(\bD^{-1/2}\bA \bD^{-1/2}) =1$, where the latter quantity denotes the spectral radius of $\bD^{-1/2}\bA \bD^{-1/2}$  and  $\|\bL^{-1}\|_{\bD}=\|\bH^{-1}\|_2$.
The inequalities in \eqref{unif-stab} then imply
\begin{align}\label{euc-2s}
\frac{1-\alpha}{1+\alpha}\,\frac{\|\br\|_{\bD}}{\|\bb\|_{\bD}} \leqslant \frac{\|\bx -\hat{\bx}\|_{\bD}}{\|\bx\|_{\bD}} 
\leqslant \frac{1+\alpha}{1-\alpha}\,\frac{\|\br\|_{\bD}}{\|\bb\|_{\bD}}
\end{align}
holds for any $0 <\alpha < 1$ and independently of the number of graph vertices where $\|\bx\|_{\bD} = \|\bD^{1/2}\bx\|_2 = \sqrt{\bx^\top \bD \bx}$ and $\br=\bb-\bL\hat{\bx}$. 

\section{Uniform Optimality} \label{unif-opt}

A consequence of the two-sided inequality \eqref{unif-stab}  is that the ratio of relative error to the relative residual is as small as $1/\kappa(\bL)< 1$ and as large as $\kappa(\bL) >1$. 
Determining the vectors $\bx$ and $\hat{\bx}$ attaining the optimal ratio is, in general, a challenging computational problem.\footnote{Note that, since $\bb=\bL\bx$, the choice of $\bx$ is equivalently characterized by the choice of $\bb$; unlike the sought $\bx$, both $\bb$ and $\hat{\bx}$ are available to compute practical \emph{a priori} and \emph{a posteriori} bounds.}

The $\bx$ and $\hat{\bx}$ attaining the optimal ratio, however, can be easily determined for the symmetric positive definite system \eqref{spd} corresponding to an undirected graph.
Decompose the solution $\bx = \bone \bgamma + \bw$ such that $\bw^\top \bD\bone = \bzero$ to conclude that $\cos_{\bD}\angle(\bx, \bone) = \|\bone \gamma\|_{\bD}/\|\bx\|_{\bD}$ and $\sin_{\bD} \angle(\bx, \bone) = \|\bw\|_{\bD}/\|\bx\|_{\bD}$. 
The equality \eqref{constant} and the decomposition $\bx = \bone \bgamma  + \bw$ imply $\bone^\top \bD \bL \bw = \bzero$, which allows us to conclude $ \| \bb \|_{\bD}^2 = \|\bone \bgamma \|_{\bD}^2 (1-\alpha)^2+ \|\bL\bw\|_{\bD}^2$.

These relationships establish:
\begin{proposition} \label{ratio-prop}
Consider the symmetric positive definite system \eqref{spd} corresponding to an undirected graph 
where $\bb$ is not collinear with $\bone$, \emph{i.e.}, $\cos_{\bD}\angle(\bb, \bone) < 1$ and the 
discount $\alpha$ satisfies $0 < \alpha < 1$. Then \eqref{euc-2s} improves to
\begin{align}\label{right-end-ptl}
\frac{1-\alpha}{1+\alpha} \frac{\|\br\|_{\bD}}{\|\bb\|_{\bD}} \leqslant \frac{\|\bx - \hat{\bx}\|_{\bD}}{\|\bx\|_{\bD}} \leqslant \kappa_{\bD}(\bL,\bb) \, \frac{\|\br\|_{\bD}}{\|\bb\|_{\bD}} \quad \bx,\bb \neq 0
\end{align}
where \[\kappa_{\bD}(\bL,\bb) \coloneqq \sqrt{\cos^2_{\bD}\angle(\bx, \bone) + \sin^2_{\bD}\angle(\bx,\bone) \big(\frac{1+\alpha}{1-\alpha}\big)^2} .\]
Moreover, if $\sin_{\bD} \angle(\bx,  \bone) \leqslant 1-\alpha$ then  $\kappa_{\bD}(L,b) \leqslant \sqrt{1+(1+\alpha)^2} \leqslant \sqrt{5}$.
\end{proposition}

The next proposition establishes computable upper bounds for the data-dependent condition number for an important class of vectors $\bb$ given by the indicator vector $\bone_{\widetilde{\Omega}}$ over a subset of vertices $\widetilde{\Omega}$ where $\bone$ is the vector of all ones and let 
$ \rho \coloneqq \bone_{\widetilde{\Omega}}^\top \bD \bone_{\widetilde{\Omega}} / \bone^\top \bD \bone $.
Let $\bL\bx = \bone - \bone_{\widetilde{\Omega}}$ and recall that solution $\bx$ is a nonnegative vector.
An elementary derivation shows that 
\begin{align*}
\cos_{\bD}(\bx, \bone) = \frac{1}{1-\alpha} \frac{\bone^\top \bD (\bone - \bone_{\widetilde{\Omega}}) }{\|\bone\|_{\bD} \| \bL^{-1} (\bone - \bone_{\widetilde{\Omega}}) \|_{\bD}} 
\geqslant \frac{\bone^\top \bD (\bone - \bone_{\widetilde{\Omega}})}{ \| \bone \|_{\bD} \| \bone - \bone_{\widetilde{\Omega}}\|_{\bD} } 
\end{align*}
since $\| \bL^{-1} (\bone - \bone_{\widetilde{\Omega}}) \|_{\bD} \leqslant \| \bL^{-1} \|_{\bD} \|\bone - \bone_{\widetilde{\Omega}}\|_{\bD} = \|\bone - \bone_{\widetilde{\Omega}}\|_{\bD}/(1-\alpha)$. 
Hence 
$\sin_{\bD}^2(\bx, \bone) \leqslant  \rho$
so that $\kappa^2_{\bD}(\bL,\bone - \bone_{\widetilde{\Omega}}) \leqslant 1 + \rho (1+\alpha)^2/(1-\alpha)^2$.  
Note that the lower bound on $\cos_{\bD}(\bx, \bone)$ is computable and is the cosine of the angle between the right-hand side 
$\bone - \bone_{\widetilde{\Omega}}$ and the vector of all ones as introduced in the abstract.
The following result summarizes our analysis:
\begin{proposition} \label{asym}
Let the fraction $\rho \coloneqq \bone_{\widetilde{\Omega}}^\top D \bone_{\widetilde{\Omega}} / \bone^\top D \bone $ and the discount $\alpha$ satisfy $0 < \alpha,\rho < 1$. 
If $\bL\bx = \bb$ with $\bb = \bone - \bone_{\widetilde{\Omega}}$  then the data-dependent condition number satisfies
\begin{align*}
\kappa_{\bD}(L,\bone - \bone_{\widetilde{\Omega}}) \leqslant 
\begin{cases}
\sqrt{5} &  \rho \leqslant (1 - \alpha)^2\\
\sqrt{1+ (1+\alpha)^2/(1-\alpha)^2}  &  \rho > (1 - \alpha)^2
\end{cases}
\end{align*}
\end{proposition}
Our conclusion is striking: when the fraction of the graph occupied by $\tilde{\Omega}$ grows modestly with the size of the graph then the relative error is essentially no larger than the relative residual error. 
When the growth is not modest, then the relative error can be arbitrarily large.
The numerical experiments we present in \S\ref{experiments} support our conclusion. 

\section{Discrete potential theory}\label{sec:dpt}
An important class of network centrality measures is given by discrete potential theory. 
The discrete potential $\bphi$ is a non-negative vector that is the solution to the linear set of equations 
\begin{align} \label{dis-pot}
\bPi^\top \bL\bPi \, \bphi  = \bc \text{ on $\Omega$ and } 
    \bphi = \bfvec\text{ on $\partial \Omega$}    
\end{align}
where  $\bc$ is a non-negative vector, $\bf $ is a boundary condition, and $\Omega \cup \partial \Omega$ is a disjoint partition of the graph. 
Let $\bPi$ contain the  columns of the order $n$ identity matrix  associated with the vertices in $\Omega$ so that $\bPi^\top \bL \bPi$ has the rows and columns of $\bL$ associated with $\Omega$. 
 
The Markov chain interpretation is given by the relationships
\begin{align}\label{prob-pot}
    \phi_i & = \begin{cases} 
                    \E_i\Big(\sum_{n=0}^T \alpha^n \, c(X_n)\Big) + \E_i\Big(f(X_T) \bone_{T< \infty}\Big) & \partial \Omega \neq \emptyset \\
                    \E_i\Big(\sum_{n=0}^\infty \alpha^n \, c(X_n) \Big)& \partial \Omega = \emptyset
                    \end{cases}
\end{align}
where $\E_i$ is expectation conditioned upon the chain starting at vertex $i$, $T$ is the random number of steps until the chain hits the boundary $\partial \Omega$. 
In words, the $i$th component of the potential averages over all discounted chains that start at vertex $i$ accumulating a cost $c_i$ when traversing vertex $i$ until the chain reaches, \emph{i.e.}, ``hits'' the boundary to incur a final cost $f$. 
The discount implies that the chain terminates with probability $1-\alpha$ per step, and the average number of steps  prior to termination is $1/(1-\alpha)$.
The vector of mean hitting-times potential can be expressed with $\alpha=1$, $\bfvec=\bzero$ and $\bc = \bPi^\top\bone_\Omega$. See \citet[\S4.2]{norr:98} for an introduction to the relationship between potential theory and Markov chains.

The bound on the inverse $\bL^{-1}$ preceding  \eqref{cond-num} then grants
$ \| \big(\bPi^\top \bL \bPi \big)^{-1} \|_p \leqslant 1/(1-\alpha)$
\citep[p.72]{govl:13}. 
The inequalities \eqref{inf-2s} then explains when the  discrete potential linear system \eqref{dis-pot} is uniformly accurate. 
Moreover, when the graph is undirected, then Proposition~\ref{ratio-prop}  applies. 
Proposition \ref{asym} applies to the important case of the mean hitting-time, $\bc = \bPi^\top\bone_\Omega$, where $\widetilde{\Omega} = \partial \Omega$ so that the relative error in the approximation $\hat{\bphi}$ is at most $\sqrt{5}$ larger than the relative residual error.

\section{PageRank} \label{pgrk}

The PageRank vector $\bz$ is given by the linear set of equations
\begin{subequations} \label{pgrk_sys}
\begin{align} \label{ppr}
  \bL^\top \bz  & = \bw (1-\alpha) 
\end{align}
where $\bA^\top \bD^{-1}$ is column stochastic and $\bw$ is a nonnegative vector.
We rewrite the equations as  $\bz^\top \bL = \bw^\top (1-\alpha)$  so that inequalities \eqref{inf-2s} hold in the $p=1$ norm because of the useful equality  $\|\bz\|_1 =\|\bz^\top\|_\infty$. The $p=1$ norm is used in PageRank analyses; see \citet[\S2]{glei:15} for a discussion.
 
We further rewrite the PageRank equations \eqref{ppr} as
\begin{align} \label{ppr-pot}
  (\bD^{-1}\bL^\top \bD) (\bD^{-1} \bz)  & = \bD^{-1} \bw (1-\alpha) \,.
\end{align}
\end{subequations}
If $\bc = \bD^{-1}\bw(1-\alpha)$ and the hitting-set $\partial \Omega$ is empty (so that $\bPi = \bI$), then by \eqref{dis-pot} the solution $\bD^{-1} \bz$ is a discrete potential on an undirected graph because  $\bD^{-1}\bL^\top \bD=\bL$. 
This grants a simple characterization of PageRank  on an undirected graph called out as missing in the literature by \citet[p.356]{glei:15}.

Let $\Omega \cup {\widetilde{\Omega}}$ denote a disjoint partition of the graph where $\Omega$ contains a small number of vertices;  as \citet{ancl:06} note, taking
$\bw = \bD\bone_\Omega (1-\alpha)^{-1} = \bD(\bone - \bone_{\widetilde{\Omega}})(1-\alpha)^{-1}$ 
induces a \emph{degree-normalized personal PageRank vector} $\bD^{-1} \bz$ which can be used for the graph cut problem. 

In this context, Proposition~\ref{ratio-prop} then bounds the relative error in an approximation to $\bw$ and Proposition~\ref{asym} then implies that the bound on the data-dependent condition number 
$\kappa_{\bD}(\bL, \bone - \bone_{\widetilde{\Omega}})$ 
is large when the fraction 
$\rho$ of the graph occupied by $\widetilde{\Omega}$ exceeds $(1-\alpha)^2$. This is typically the case because $\Omega$ is chosen to be a small set and hence the the complementary fraction $1-\rho$ of the graph occupied by $\Omega$ is also small, implying $\rho$ is close to $1 \gg (1-\alpha)^2$.

The next result explains that the ensuing ranking of the vertices is unaffected whether 
$\bb=\bone_{\widetilde{\Omega}}$ or $\bb=\bone-\bone_{\widetilde{\Omega}}$.
\begin{proposition}\label{orderings}
The orderings induced by the solutions to the two linear systems $\bL\bx = \bone - \bone_{\widetilde{\Omega}}$ and 
$\bL\bu = \bone_{\widetilde{\Omega}}$ are the same.
\end{proposition}
The result follows immediately because $\bone = \bL(\bx+\bu) $ so that \eqref{constant} implies $\bx + \bu = \bone/(1+\alpha)$. 
Hence, the orderings, \emph{i.e.}, the indices when the elements of $\bx$ are sorted in ascending order are the same as the indices that arise when the elements of $\bu$ are sorted in descending order because both $\bx$ and $\bu$ are nonnegative vectors. 
The practical conclusion is immediate: Solve $\bL\bx = \bone - \bone_{\widetilde{\Omega}}$ when $\rho \leqslant (1-\alpha)^2$ and otherwise solve $\bL\bu = \bone_{\widetilde{\Omega}}$  to obtain an accurate approximation inducing a reliable ordering of vertices by centrality.

\begin{figure*}[htb]
\begin{centering}
\begin{tabular}{cc}
\includegraphics[width=\textwidth]{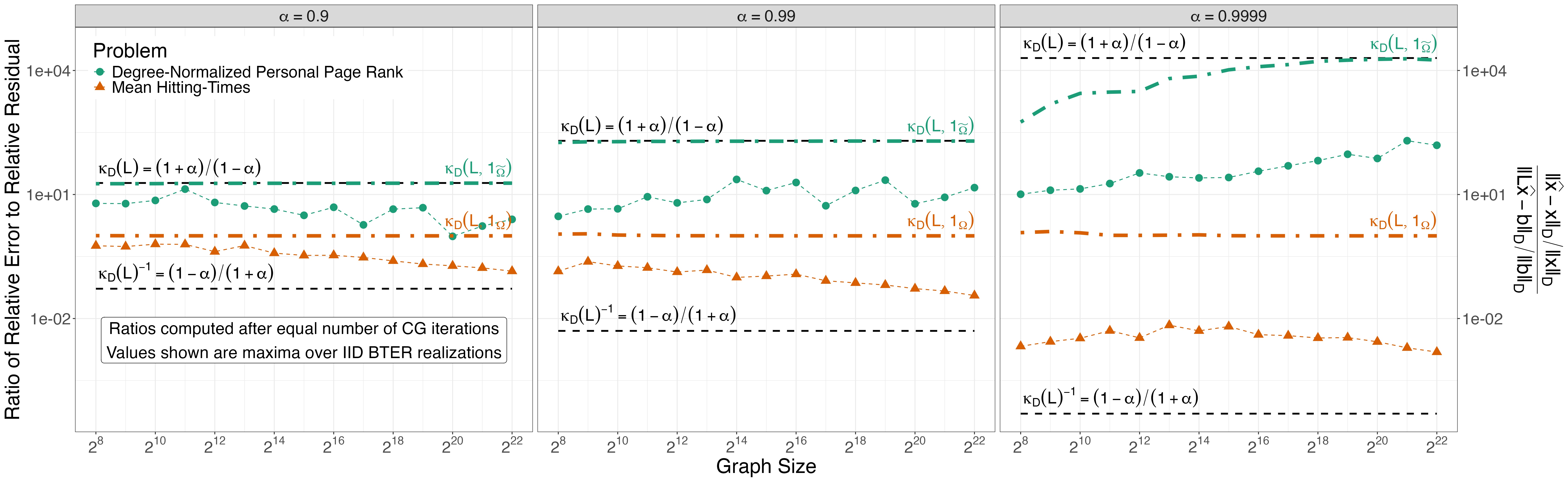}   
\end{tabular}
\caption{\label{fig:G2-G3} 
Three plots of the ratio of relative error to the norm of the relative residual over three choices of discount of $\alpha$ on a sample of BTER undirected graphs with the choice of $\bb=\bone_{\Omega}$ (MHT) and $ b= \bone_{\widetilde{\Omega}}$ (PPR) 
where $\widetilde{\Omega}$ contains a vertex of the graph selected randomly and $\bone_{\widetilde{\Omega}}+\bone_{\Omega} = \bone$.  The solution $x$ was computed via the conjugate gradient (CG) iteration using a relative residual tolerance of $10^{-12}$. The approximations $\hat{\bx}$ were computed via the CG iteration terminated when the relative residual tolerance satisfied $10^{-3}$ was achieved or a maximum of $40$, $120$, and $480$ CG iterations were performed. 
}
\end{centering}
\end{figure*}

\section{Experiments}\label{experiments}

Figure~\ref{fig:G2-G3} depicts the dramatic impact the choices of $\bb=\bone - \bone_{\widetilde{\Omega}}$ and $\bb=\bone_{\widetilde{\Omega}}$ implied by Propositions~\ref{ratio-prop}--\ref{asym} has for three different choices of discount $\alpha$ when $\widetilde{\Omega}$ contains a vertex of the graph selected randomly.
We computed approximations to the
mean hitting-time (MHT) and the degree-normalized personal PageRank (PPR) vectors on the random undirected
graph model BTER \citep{bter} with heavy-tailed degree and clustering coefficient distributions using the conjugate
gradient iteration \citep{cg-iter}. 
We report the maximum observed ratio of the relative error to the relative residual error, 
which is an estimate of $\kappa_{\bD}(\bL,\bb)$, in order to highlight the effects of discounting and of the choice of $\bb$ on the conditioning of each problem.
 
The behavior of the PPR problem (green circles) correlates well with the condition number bound $(1+\alpha)/(1-\alpha)$ (black line) in stark contrast to the MHT problem (orange triangles). 
The data-dependent condition numbers (green and orange dashed lines) underscore the improvement concluded by Proposition \ref{asym} for the MHT problem---for graphs of increasing size, the relative error is guaranteed to be no larger than $\sqrt{5}$ times the relative residual error, a useful certificate of solution quality.  

In particular, because $\bb = \bone_{\widetilde{\Omega}}$ is nearly orthogonal to $\bone$, the PPR  data-dependent
condition number provides only marginal improvement over the general bound in Proposition~\ref{ratio-prop}, with the disparity between MHT and PPR increasing rapidly with $\alpha$. This disparity is remarkable when we recall that Proposition \ref{orderings} implies the two systems extract the same ranking of vertices. Finally, we emphasize that these trends depend upon the properties of the linear set of equations and hold regardless of the choice of algorithm used to solve that system.
 
\section*{Acknowledgements}
We acknowledge helpful discussions with Daniel M. Dunlavy, Mike Eydenberg, Renee Gooding, Alex Foss, J. Derek Tucker and Natalie Wellen.

Sandia National Laboratories is a multimission laboratory managed and operated by National Technology \&
Engineering Solutions of Sandia, LLC, a wholly owned subsidiary of Honeywell International Inc., for the U.S.
Department of Energy’s National Nuclear Security Administration under contract DE-NA0003525. This paper describes objective technical results and analysis. Any subjective views or opinions that might be expressed in the paper do not necessarily represent the views of the U.S. Department of Energy or the United States Government. 

Approved for Release: SAND2024-05956O

\printbibliography

\end{refsection}
\end{document}